\documentclass[12pt]{article}

\usepackage{latexsym}
\usepackage{amssymb}

\newcommand{\ga}{\alpha}
\newcommand{\gb}{\beta}
\renewcommand{\gg}{\gamma}
\newcommand{\gd}{\delta}

\newcommand{\gk}{\kappa}
\newcommand{\gl}{\lambda}

\newcommand{\gs}{\sigma}

\newcommand{\go}{\omega}

\newcommand{\la}{\langle}
\newcommand{\ra}{\rangle}
\newcommand{\ov}{\overline}

\newcommand{\ha}{\aleph}

\newcommand{\forces}{\Vdash}

\renewcommand{\models}{\vDash}

\newcommand{\dom}{{\rm dom}}

\newcommand{\FC}{{\mathbb C}}

\newcommand{\FP}{{\mathbb P}}
\newcommand{\FQ}{{\mathbb Q}}
\newcommand{\FR}{{\mathbb R}}
\newcommand{\FS}{{\mathbb S}}
\newcommand{\FT}{{\mathbb T}}
\newcommand{\K}{{\cal K}}

\newtheorem{theorem}{Theorem}

\newtheorem{prop}[theorem]{Proposition}

\newtheorem{lemma}[theorem]{Lemma}

\newtheorem{Main Question}[theorem]{Main Question}
\newenvironment{proof}{\noindent{\bf
Proof:}}{\nopagebreak\mbox{}\newline\makebox[\textwidth]{\hfill$\square$}
\par\bigskip}

\def\P{{\mathbb P}}

\def\intersect{\cap}

\def\and{\mathrel{\kern1pt\&\kern1pt}}
\def\image{\mathbin{\hbox{\tt\char'42}}}

\def\<#1>{\langle\,#1\,\rangle}

\newcommand{\cp}{\mathop{\rm cp}\nolimits}
\title{Exactly Controlling the Non-Supercompact Strongly Compact
Cardinals\thanks{
Supercompact cardinal, strongly compact cardinal,
strong cardinal, indestructibility,
non-reflecting stationary set of ordinals}}
\author{Arthur W.~Apter\thanks{The first author's research was partially
        supported by PSC-CUNY Grants 61449-00-30 and
        63436-00-32.}\\
        Department of Mathematics\\
        Baruch College of CUNY\\
        New York, New York 10010 USA\\
        http://math.baruch.cuny.edu/$\sim$apter\\
        awabb@cunyvm.cuny.edu\\
        \\
        Joel David Hamkins\thanks{The second author's research was
        partially supported
        by PSC-CUNY and NSF grants.}\\
        The College of Staten Island of CUNY and\\
        The CUNY Graduate Center\\
        http://jdh.hamkins.org\\
        jdh@hamkins.org}

\date{February 25, 2002\\
      (revised December 30, 2002)}

\begin{document}
\maketitle

\begin{abstract}We summarize the known methods of producing a non-supercompact
strongly compact cardinal and describe some new
variants. Our Main Theorem shows how to apply
these methods to many cardinals simultaneously
and exactly control which cardinals are
supercompact and which are only strongly compact
in a forcing extension. Depending upon the
method, the surviving non-supercompact strongly
compact cardinals can be strong cardinals, have
trivial Mitchell rank or even contain a club
disjoint from the set of measurable cardinals.
These results improve and unify Theorems 1 and 2
of \cite{A97}, due to the first author.
\end{abstract}
\baselineskip=24pt

\section{Introducing the Main Question}\label{s1}

The notions of strongly compact and supercompact
cardinal are very close, so close that years ago
it was an open question whether they were
equivalent. When Solovay first defined the
supercompact cardinals---witnessed by {\it
normal} fine measures on $P_\kappa(\gl)$, rather
than merely fine measures, which give rise only
to strong compactness---he had simply added the
kind of normality assumption that set theorists
were accustomed to getting for free in the case
of measurable cardinals. And so he and many
others expected that these notions would be
equivalent. Confirming evidence for this initial
expectation was found in the models constructed
by Kimchi and Magidor \cite{KM} (see \cite{A98}
and \cite{AS97a} for a modern account), where the
two notions do coincide.

But set theorists now know that the two notions
are not equivalent (although it is unknown
whether they are equiconsistent), and we have a
variety of ways of producing non-supercompact
strongly compact cardinals. Let us classify these
different methods into five general categories.

\begin{enumerate}

\item {\it The purely ZFC method.} This is historically the
first way in which non-supercompact strongly
compact cardinals were exhibited, and can be
contrasted with the forcing results, which show
merely the relative consistency of inequivalence.
Its foundation is the observation by Solovay's
student Telis Menas \cite{Me} that any measurable
limit of strongly compact cardinals is itself
strongly compact, but the least such cardinal
(and many more) cannot be supercompact. A related
observation, given in the context of a
supercompact limit of supercompact cardinals by
the first author in \cite{A80}, is that the least
cardinal $\gk$ (and many more) that is $\gk^+$,
$\gk^{++}$, $\gk^{+++}$, etc$.$ supercompact and
is also a limit of strongly compact cardinals is
strongly compact but isn't (fully)
supercompact.\footnote{The strong compactness
measures for the limit cardinal $\kappa$ are
obtained simply by integrating the smaller strong
compactness measures for cardinals below $\kappa$
with respect to a fixed measure on $\kappa$.
Conversely, the least such cardinal $\kappa$ must
have limited supercompactness in $M$ for any
embedding $j:V\to M$ having critical point
$\kappa$, by the minimality of $j(\kappa)$ there;
so $\kappa$ cannot be supercompact.}

\item {\it The Menas forcing method.}
Menas used his aforementioned result in \cite{Me}
to exhibit the inequivalence of strong
compactness and supercompactness much lower in
the hierarchy, by forcing over a model with a
measurable limit of supercompact cardinals and
producing a model where the least strongly
compact cardinal is not supercompact. Jacques
Stern, in unpublished work, generalized these
ideas to create a model in which the first two
strongly compact cardinals are not supercompact.
The first author, in \cite{A97}, also used these
ideas in a way we will discuss in greater detail
below. The second author combined these ideas
with fast function forcing in Corollary 4.3 of
\cite{H4}, where he showed that {\it any}
strongly compact cardinal $\kappa$ can be forced
to be a non-supercompact strongly compact
cardinal. This was accomplished by showing that
any strongly compact cardinal $\kappa$ can be
made indestructible by the forcing to add a club
$C\subseteq\kappa$ containing no measurable
cardinals. After such forcing, $\kappa$ clearly
cannot be supercompact, or even have nontrivial
Mitchell rank. The technique allows one to add
coherent sequences of clubs to smaller cardinals
which reflect at their inaccessible limit points.

\item {\it The method of iterated Prikry forcing.}
Inspired by Menas' work mentioned in (1) and (2)
above, Magidor \cite{Ma} provided a technique for
producing non-supercompact strongly compact
cardinals, the method of iterated Prikry forcing,
which yielded a striking result: a model in which
the least measurable cardinal---which obviously
cannot be supercompact---is strongly compact.
Magidor's ground model required no GCH
assumptions and began with only a strongly
compact cardinal. These ideas have also been used
by the first author in \cite{A80} and \cite{A97a}
and by Abe in \cite{Ab} to create further
examples of non-supercompact strongly compact
cardinals. In addition, a modification of this
forcing notion was discovered by Gitik in
\cite{G}, where a different proof of Magidor's
theorem is essentially given in Section 4, pages
302--303, starting with a ground model satisfying
GCH and containing a strongly compact cardinal,
using an iteration of Prikry forcing containing
Easton supports (as opposed to finite support in
the first coordinate and full support in the
second coordinate in the Magidor iteration of
Prikry forcing). Gitik's technique was later
modified by the first author and Gitik in
\cite{AG}, where, starting from a ground model
satisfying GCH and containing a supercompact
cardinal, they produced an assortment of models
in which the least strongly compact cardinal
$\gk$ isn't supercompact yet has its strong
compactness and degree of supercompactness fully
indestructible under $\gk$-directed closed
forcing. In one of the models constructed, the
least strongly compact cardinal is also the least
measurable cardinal, yet indestructible.

\item {\it The method of iteratively adding non-reflecting stationary
sets of ordinals.} While Magidor's result
mentioned in (3) above was very exciting, making
strong compactness and supercompactness seem very
far apart, it unfortunately could not be used to
handle more than one strongly compact cardinal.
This was because Prikry forcing above a strongly
compact cardinal adds a weak square sequence,
which destroys the strong compactness of the
smaller cardinal. Magidor overcame this
difficulty by inventing yet another technique for
producing non-supercompact strongly compact
cardinals. Rather than iterating Prikry forcing
below a supercompact cardinal $\kappa$, he now
iterated instead the forcing to add a
non-reflecting stationary set of ordinals to
every measurable cardinal below $\kappa$. This
iteration destroys all the measurable cardinals
below $\kappa$, and yet $\kappa$ remains strongly
compact in the extension. Then, using this
technique, Magidor constructed for each natural
number $n$ a model where the first $n$ measurable
cardinals are strongly compact (but clearly not
supercompact), beginning with a model having $n$
indestructible supercompact cardinals.
Unfortunately, Magidor's method does not seem to
work in the infinite case, and the question of
the relative consistency of the first $\go$
measurable cardinals all being strongly compact
is open. This theorem and technique, although
unpublished by Magidor, appeared in \cite{AC1},
along with a related generalization due to the
first author and Cummings. Further
generalizations of Magidor's method using an
iteration of the forcing for adding
non-reflecting stationary sets of ordinals to
produce additional examples of non-supercompact
strongly compact cardinals can be found in
\cite{A95} and \cite{AC2}, and a modification of
this method using an iteration of the forcing for
adding Cohen subsets to produce additional
examples of non-supercompact strongly compact
cardinals can be found in \cite{AH4}.
Another exposition of Magidor's method can be
found in
\cite{A01}, as well as Lemma \ref{l1} of this
paper.

\item {\it The method of Radin forcing.}
This unpublished technique is due to Woodin, and
was invented by him and modified by Magidor at
the January 7-13, 1996 meeting in Set Theory held
at the Mathematics Research Institute,
Oberwolfach, Germany to give a proof of the
theorem of \cite{A97a} from only one supercompact
cardinal, instead of hypotheses on the order of a
supercompact limit of supercompact cardinals. A
further proof of this theorem, employing only one
supercompact cardinal, is given in \cite{AC1}.

\end{enumerate}

We mention all these methods for turning a
supercompact cardinal into a non-supercompact
strongly compact cardinal because in this article
we seek uniform methods to control exactly in
this way any given class of supercompact
cardinals simultaneously. Specifically, we seek
an affirmative answer to the following question.

\begin{Main Question}\label{MainQuestion}

Given a class of supercompact cardinals, can one
force them to be strongly compact and not
supercompact while fully preserving all other
supercompact cardinals and creating no new
strongly compact or supercompact cardinals?

\end{Main Question}

We are pleased to announce an affirmative answer
to this question for a broad collection of
classes. Yes, one can force any given class of
supercompact cardinals $\cal A$ to become
strongly compact and not supercompact while fully
preserving all other supercompact cardinals and
creating no new strongly compact or supercompact
cardinals, provided that $\cal A$ does not
contain certain kinds of complicated limit
points. In particular, if there are no
supercompact limits of supercompact cardinals (or
merely no supercompact limits of supercompact
limits of supercompact cardinals), then any class
$\cal A$ of supercompact cardinals can be exactly
controlled in this way. Indeed, our control over
the non-supercompact strongly compact cardinals
is even greater than requested in the Main
Question, for we can ensure that they remain
strong cardinals in the extension or
alternatively, that as measurable cardinals they
have trivial Mitchell rank or even contain a club
containing no measurable cardinals. Specifically,
we will prove the following.

\newtheorem{MainTheorem}[theorem]{Main Theorem}
\begin{MainTheorem}\label{t1}

Suppose that ${\cal A}$ is a subclass of the
class $\K$ of supercompact cardinals containing
none of its limit points. Then there is a forcing
extension $V^\FP$ in which the cardinals in
${\cal A}$ remain strongly compact but become
non-supercompact, while the cardinals in $\K -
{\cal A}$ remain fully supercompact. In addition,
no new strongly compact or supercompact cardinals
are created. In $V^\FP$, the class of strongly
compact cardinals is composed of $\K$ together
with its measurable limit points. Depending on
the choice of $\FP$, the cardinals of ${\cal A}$
become strong cardinals in $V^\FP$ or contain a
club disjoint from the measurable cardinals
there, respectively. Finally, if $\gk$ is a
measurable limit of ${\cal A}$ at which the GCH
holds in $V$, then both the GCH at $\gk$ and
$\gk$'s measurability are preserved in $V^\FP$.

\end{MainTheorem}

\newtheorem{MainCorollary}[theorem]{Main Corollary}
\begin{MainCorollary}

If there is no supercompact limit of supercompact
cardinals, then the answer to the Main Question
is Yes. In particular, an affirmative answer to
the Main Question is relatively consistent with
the existence of many supercompact cardinals,
even a proper class of supercompact cardinals.

\end{MainCorollary}

\noindent The point here is that if there is no
supercompact limit of supercompact cardinals,
then every class $\cal A$ satisfies the
hypothesis of the theorem.

By iterating the result of the Main Theorem, we
are able to generalize it to make almost the same
conclusions with any class $\cal A$ having finite
Cantor-Bendixon rank.

\newtheorem{Generalization}[theorem]{Generalized Main Theorem}
\begin{Generalization}

The main conclusions of the Main Theorem \ref{t1}
hold for any class ${\cal A}\subseteq{\cal K}$
having finite Cantor-Bendixon rank. That is, for
any such ${\cal A}$, there is a forcing extension
$V^\FP$ in which the cardinals of ${\cal A}$
become strongly compact but not supercompact,
while the cardinals in $\K - {\cal A}$ remain
fully supercompact. In addition, no new strongly
compact or supercompact cardinals are created. In
$V^\FP$, the class of strongly compact cardinals
is composed of $\K$ together with its measurable
limit points. Finally, if $\gk$ is a measurable
limit of ${\cal A}$ at which the GCH holds in
$V$, then both the GCH at $\gk$ and $\gk$'s
measurability are preserved in $V^\FP$.

\end{Generalization}

And the Generalized Main Theorem comes also with
its own generalized corollary:

\newtheorem{GeneralizedMainCorollary}[theorem]{Generalized Main Corollary}
\begin{GeneralizedMainCorollary}

If the class of supercompact cardinals has finite
Cantor-Bendixon rank, then the answer to the Main
Question \ref{MainQuestion} is Yes. Therefore, an
affirmative answer to the question is relatively
consistent with the existence of a proper class
of supercompact limits of supercompact limits of
supercompact cardinals, and supercompact limits
of these, and so on.

\end{GeneralizedMainCorollary}

The first author made substantial progress in
\cite{A97} on a question related to the Main
Question, when he proved that we already have a
fine control over the pattern that the
supercompact cardinals make as a subclass of the
strongly compact cardinals. Specifically, in
Theorem 1 of \cite{A97}, using ideas of Menas
from \cite{Me} in tandem with those from
\cite{A98}, he begins with an inaccessible limit
$\Omega$ of measurable limits of supercompact
cardinals, and produces a model where the pattern
of the supercompact cardinals as a subset of the
general class of (strongly) compact cardinals
follows any prescribed function $f:\Omega\to 2$
in the ground model.\footnote{The proof could be
modified to omit the cardinal $\Omega$, and treat
$f:Ord\to 2$, through the use of proper classes.}
The resulting non-supercompact strongly compact
cardinals there have trivial Mitchell rank. In
his argument, somewhat like the earlier arguments
of Menas and Stern, the cardinals that eventually
become strongly compact but not supercompact
begin as measurable limits of supercompact
cardinals in the ground model. Necessarily,
therefore, many supercompact cardinals are
destroyed along the way. The main result of this
article overcomes this difficulty.

Our Main Theorem provides a new proof of Theorem
1 of \cite{A97} and provides a substantial
reduction in the hypotheses used to prove Theorem
1 of \cite{A97}, from an inaccessible limit of
measurable limits of supercompact cardinals to a
proper class of supercompact cardinals. It also
allows us to provide a uniform proof of both a
strengthened version of Theorem 1 of \cite{A97}
and a more general, stronger version of Theorem 2
of \cite{A97}, since Theorem 2 of \cite{A97}
extends Theorem 1 of \cite{A97} to the situation
encompassing a cardinal $\gk$ which is a
supercompact limit of supercompact cardinals.
Further, Theorem 1 of \cite{A01} can now be
derived as a corollary of our Main Theorem.

Let us conclude this section with some
preliminary information and basic definitions.
Essentially, our notation and terminology are
standard, and when this is not the case, it will
be clearly noted. For $\ga < \gb$ ordinals,
$[\ga, \gb], [\ga, \gb), (\ga, \gb]$, and $(\ga,
\gb)$ are as in standard interval notation.

When forcing, $q \ge p$ will mean that $q$ is
stronger than $p$. If $G$ is $V$-generic over
$\FP$, we will use both $V[G]$ and $V^{\FP}$ to
indicate the universe obtained by forcing with
$\FP$. If $\FP$ is an iteration, then $\FP_\ga$
is the forcing up to stage $\ga$. When $\gk$ is
inaccessible and $\FP = \la \la \FP_\ga, \dot
\FQ_\ga \ra : \ga < \gk \ra$ is an Easton support
iteration of length $\gk$ which at stage $\ga$
performs some nontrivial forcing based on the
ordinal $\gd_\ga$, then we will say that
$\gd_\ga$ is in the field of $\FP$. If $x \in
V[G]$, then $\dot x$ will be a term in $V$ for
$x$. We may, from time to time, confuse terms
with the sets they denote and write $x$ when we
actually mean $\dot x$, especially when $x$ is
some variant of the generic set $G$, or $x$ is in
the ground model $V$.

If $\gk$ is a cardinal and $\FP$ is a partial
ordering, $\FP$ is $\gk$-directed closed if for
every directed subset $D\subseteq\FP$ of size
less than $\kappa$ (where $D$ is directed if
every two elements of $D$ have an upper bound in
$\FP$) has an upper bound in $\FP$. The partial
order $\FP$ is $\gk$-strategically closed if in
the two person game in which the players
construct an increasing sequence $\langle p_\ga:
\ga \le\gk\rangle$, where player I plays odd
stages and player II plays even and limit stages
(choosing the trivial condition at stage 0), then
player II has a strategy which ensures the game
can always be continued.
Note that if $\FP$ is $\gk^+$-directed closed,
then $\FP$ is $\gk$-strategically closed. Also,
if $\FP$ is $\gk$-strategically closed and $f :
\gk \to V$ is a function in $V^\FP$, then $f \in
V$. $\FP$ is ${\prec}\gk$-strategically closed if
in the two person game in which the players
construct an increasing sequence $\langle p_\ga:
\ga < \gk\rangle$, where player I plays odd
stages and player II plays even and limit stages
(again choosing the trivial condition at stage
0), then player II has a strategy which ensures
the game can always be continued.

If $X$ is a set of ordinals, then $X'$ is the set
of limit points of $X$. The Cantor-Bendixon
derivatives of a set $X$ are defined by
iteratively removing isolated points. One begins
with the original set $X^{(0)}=X$, removes
isolated points at successor stages by keeping
only the limit points
$X^{(\alpha+1)}=X^{(\alpha)}\intersect
(X^{(\alpha)})'$, and takes intersections at
limit stages
$X^{(\lambda)}=\bigcap_{\alpha<\lambda}
X^{(\alpha)}$. For any class of ordinals $X$, if
$X^{(\alpha)}=\emptyset$ for some $\alpha$, then
the Cantor-Bendixon rank of $X$ is the least such
$\alpha$. The Cantor-Bendixon rank of a point
$\gamma$ in $X$ is the largest $\beta$ such that
$\gamma\in X^{(\beta)}$, that is, the stage at
which $\gamma$ becomes isolated.

In this paper, we will use non-reflecting
stationary set forcing $\FP_{\eta,\gl}$.
Specifically, if $\eta<\lambda$ are both regular
cardinals, then conditions in $\FP_{\eta,\gl}$
are bounded subsets $s\subset\gl$ consisting of
ordinals of cofinality $\eta$ such that for every
$\alpha<\lambda$, the initial segment
$s\intersect\alpha$ is non-stationary in
$\alpha$, ordered by end-extension. It is
well-known that if $G$ is $V$-generic over
$\FP_{\eta, \gl}$ (see \cite{Bu}, \cite{A01}, or
\cite{KM}) and the GCH holds in $V$, then in
$V[G]$, the set $S=S[G]=\bigcup G \subseteq \gl$
is a non-reflecting stationary set of ordinals of
cofinality $\eta$, the bounded subsets of $\gl$
are the same as those in $V$, and cardinals,
cofinalities and the GCH have been preserved. It
is virtually immediate that $\FP_{\eta, \gl}$ is
$\eta$-directed closed. It follows from Theorem
4.8 of \cite{SRK} that the existence of a
non-reflecting stationary subset of $\lambda$,
consisting of ordinals of confinality $\eta$,
implies that no cardinal
$\delta\in(\eta,\lambda]$ is $\lambda$ strongly
compact. Thus, iterations of this forcing provide
a way to destroy all strongly compact cardinals
in an interval.

We assume familiarity with the large cardinal
notions of measurability, strongness, strong
compactness, and supercompactness. Interested
readers may consult \cite{K} for further details.
We mention only that a cardinal $\gk$ is ${<}\gl$
supercompact iff it is $\gd$ supercompact for
every cardinal $\gd < \gl$. We will always
identify an ultrapower with its Mostowski
collapse. We note that a measurable cardinal
$\gk$ has trivial Mitchell rank if there is no
embedding $j : V \to M$ for which $\cp(j)=\kappa$
and $M \models ``\gk$ is measurable''. An
ultrafilter ${\cal U}$ generating this sort of
embedding will be said to have trivial Mitchell
rank as well. Ultrafilters of trivial Mitchell
rank always exist for any measurable cardinal
$\gk$. Also, unlike \cite{K}, we will say that
the cardinal $\gk$ is $\gl$ strong for an ordinal
$\gl > \gk$ if there is $j : V \to M$ an
elementary embedding having critical point $\gk$
so that $j(\gk)
> |V_\gl|$ and $V_\gl \subseteq M$. As always, $\gk$ is strong if $\gk$
is $\gl$ strong for every $\gl > \gk$.

As in \cite{H4} we define the lottery sum of a
collection $\cal C$ of partial orderings to be
$\oplus{\cal C}=\{\la \FQ, q \ra : \FQ \in {\cal
C}$ and $q \in \FQ\} \cup \{0\}$, ordered with
$0$ below everything and $\la \FQ, q \ra \le \la
\FQ', q' \ra$ iff $\FQ = \FQ'$ and $q \le q'$.
(This is equivalent simply to taking the product
of the corresponding Boolean algebras.) Forcing
with $\oplus\cal C$ amounts to selecting a
particular partial ordering $\FQ$ from $\cal C$,
the ``winning partial ordering'', and then
forcing with it. The lottery preparation of
\cite{H4} proceeds by iterating these lottery
sums, and has proved to be useful for obtaining
indestructibility even in a large cardinal
context in which one lacks a Laver function.

Finally, we will say that an elementary embedding
$k : V \to N$ with critical point $\gk$ has the
$\lambda$-cover property when for any $x
\subseteq N$ with $|x| \le \gl$, there is some $y
\in N$ so that $x \subseteq y$ and $N \models
``|y| < k(\gk)$''. A suitable cover of
$j\image\lambda$ generates a fine measure over
$P_\kappa(\lambda)$ and conversely, so one can
easily deduce that such an embedding exists iff
$\kappa$ is $\lambda$ strongly compact (see
Theorem 22.17 of \cite{K}).

\section{Two Useful Propositions}\label{s2}

The proof of the Main Theorem will proceed by
taking large products of the forcing that
transforms a given supercompact cardinal into a
non-supercompact strongly compact cardinal. The
particular iterations we will use to accomplish
this are given in the proofs of the following two
propositions.

\begin{prop}\label{p1}

If $\gk$ is supercompact, then regardless of the
number of large cardinals in the universe, there
is a forcing extension $V^\P$ in which $\gk$ is
strongly compact, $2^\gk = \gk^+$ and $\gk$ has
trivial Mitchell rank. For any regular $\eta <
\gk$, such a partial ordering $\FP$ can be found
which is $\eta$-directed closed and of
cardinality at most $2^\gk$. Indeed, if $2^\gk =
\gk^+$ in the ground model, then $\FP$ can have
cardinality $\gk$. Furthermore, $\FP$ can be
defined so as to destroy any strongly compact
cardinal in the interval $(\eta, \gk)$.

\end{prop}

\begin{prop}\label{p2}

Suppose that $\gk$ is supercompact and $\eta <
\gk$. Then regardless of the number of large
cardinals in the universe, there is a forcing
extension $V^\FP$ in which $\kappa$ becomes a
non-supercompact strongly compact cardinal,
$2^\kappa=\kappa^+$ and all other supercompact
cardinals above $\eta$ are preserved. Such a
partial ordering $\FP$ can be found which is
$\eta$-directed closed and of cardinality at most
$2^\gk$. Indeed, if $2^\gk = \gk^+$ in the ground
model, then $\FP$ can have cardinality $\gk$.
Depending upon the exact choice of $\FP$, the
cardinal $\gk$ will either contain a club
disjoint from the measurable cardinals (and hence
have trivial Mitchell rank) or become a strong
cardinal, respectively. Furthermore, every
strongly compact cardinal in $V^{\FP}$ in the
interval $(\eta, \gk]$ is either supercompact in
$V$ or a measurable limit of supercompact
cardinals in $V$.

\end{prop}

These two propositions are closely related to
\cite[Corollary 4.3]{H4}. In particular, if one
omits the last sentence of Proposition \ref{p1}
(which will not actually be relevant in our
application), it is an immediate consequence of
\cite[Corollary 4.3]{H4}, which has both a weaker
hypothesis and a stronger conclusion: one can add
to any strongly compact cardinal $\kappa$ a club
disjoint from the measurable cardinals while
preserving the strong compactness of $\kappa$ and
neither creating nor destroying any measurable
cardinals. This method also arises in the proof
of Proposition \ref{p2}, in Lemma \ref{n1}. And
the forcing of \cite[Corollary 4.3]{H4} can be a
component of the product forcing used to prove
the Main Theorem in the situation when there are
no supercompact limits of supercompact cardinals
(which is also true for the forcing given in
Proposition \ref{p1} or the forcing used in the
proof of Theorem 1 of \cite{AC2}). We give the
alternative proof of Proposition \ref{p1} here
because of the extra property that it destroys
all strongly compact cardinals in the interval
$(\eta,\kappa)$, which may find an application
elsewhere.

We would like to call special attention to the
fact that both Propositions \ref{p1} and \ref{p2}
can be applied over a universe with many large
cardinals. This contrasts sharply with Magidor's
forcing to create a non-supercompact strongly
compact cardinal by adding non-reflecting
stationary sets of ordinals to every measurable
cardinal below the supercompact cardinal $\gk$ or
the modifications of this method given in
\cite{A95} and \cite{AC1}, which seem to require
severe restrictions on the type of large
cardinals above $\gk$.\footnote{The modification
of Magidor's method given in \cite{AC2} also can
be applied over a universe with many large
cardinals.} Also, the forcing of Proposition
\ref{p2} has been explicitly designed to preserve
the supercompactness of all cardinals above
$\eta$ except for $\gk$, which distinguishes it
from the other partial orderings we have
mentioned.

Let's now prove Proposition \ref{p1}.

\begin{proof}
Let $V \models ``$ZFC + $\gk$ is supercompact'',
and suppose $\eta < \gk$ is regular. By forcing
if necessary, we may assume without loss of
generality that $2^\gk = \gk^+$ in $V$. This is
done by simply using the Laver preparation from
$\eta$ up to $\kappa$ (e.g., the version given in
\cite{A98}), followed by the forcing which adds a
Cohen subset to $\gk^+$, thereby ensuring
$2^\kappa=\kappa^+$. This combined forcing also
preserves all supercompact cardinals above $\eta$
and has cardinality $2^\gk$. By the Gap Forcing
Theorem of \cite{H2} and \cite{H3}, if this
combined forcing begins by adding a Cohen subset
to $\eta$, it creates no new supercompact or
measurable cardinals above $\eta$.

Let $\FP$ be the Easton support forcing
$\kappa$-iteration which adds to every measurable
limit of strong cardinals $\gs \in (\eta, \gk)$ a
non-reflecting stationary set of ordinals of
cofinality $\eta$. It is not difficult to see
that $\FP$ has cardinality $\gk$, so $V^\FP
\models ``2^\gk = \gk^+$''. Thus, the following
three lemmas complete the proof of Proposition
\ref{p1}.

\begin{lemma}\label{l1}

$V^\FP \models ``\gk$ is strongly compact''.

\end{lemma}

\begin{proof}
We use Magidor's method for preserving strong
compactness mentioned at the beginning of this
paper. Let $\gl > \gk$ be an arbitrary singular
strong limit cardinal of cofinality at least
$\gk$, and let $k_1 : V \to M$ be an elementary
embedding witnessing the $\gl$ supercompactness
of $\gk$ so that $M \models ``\gk$ isn't $\gl$
supercompact''. By the choice of $\gl$, the
cardinal $\kappa$ is ${<}\lambda$ supercompact in
$M$. Since $\lambda\geq 2^\kappa$, we know
$\kappa$ is measurable in $M$. Therefore, there
is a normal measure of trivial Mitchell rank over
$\kappa$ in $M$, yielding an embedding $k_2 : M
\to N$, with critical point $\kappa$, such that
$N \models ``\gk$ isn't measurable''. In
addition, as $\gl \ge 2^\gk$, Lemma 2.1 of
\cite{AC2} and the succeeding remark imply that
in both $V$ and $M$, $\gk$ is a strong cardinal
which is also a limit of strong cardinals, and in
fact, in both $V$ and $M$, $\gk$ carries a normal
measure concentrating on strong cardinals.

It is easy to verify that the composed embedding
$j = k_2 \circ k_1:V\to N$ has the
$\lambda$-cover property, and therefore witnesses
the $\gl$ strong compactness of $\gk$. We will
show that $j$ lifts to $j : V^{\FP} \to
N^{j(\FP)}$. Since this lifted embedding will
witness the $\gl$ strong compactness of $\gk$ in
$V^{\FP}$, this will prove Lemma \ref{l1}.

To do this, factor $j(\FP)$ as $\FP \ast \dot \FQ
\ast \dot \FR$, where $\dot \FQ$ is a term for
the portion of $j(\FP)$ from stage $\kappa$ up to
and including stage $k_2(\gk)$, and $\dot \FR$ is
a term for the rest of $j(\FP)$, from stage
$k_2(\kappa)+1$ up to $j(\kappa)$. Since $N
\models ``\gk$ isn't measurable'', we know that
$\gk \not\in {\rm field}(\dot \FQ)$. Thus, the
field of $\dot \FQ$ is composed of all
$N$-measurable limits of $N$-strong cardinals in
the interval $(\gk, k_2(\gk)]$ (so $k_2(\gk)$ is
in the field of $\dot \FQ$), and the field of
$\dot \FR$ is composed of all $N$-measurable
limits of $N$-strong cardinals in the interval
$\bigl(k_2(\gk), k_2(k_1(\gk))\bigr)$.

Let $G_0$ be $V$-generic over $\FP$. We will
construct in $V[G_0]$ an $N[G_0]$-generic object
$G_1$ over $\FQ$ and an $N[G_0][G_1]$-generic
object $G_2$ over $\FR$. Since $\FP$ is an Easton
support iteration of small forcing, with a direct
limit at stage $\kappa$ and no forcing right at
stage $\gk$, the construction of $G_1$ and $G_2$
ensures that $j '' G_0 \subseteq G_0 \ast G_1
\ast G_2$. It follows that $j : V \to N$ lifts to
$j : V[G_0] \to N[G_0][G_1][G_2]$ in $V[G_0]$.

To build $G_1$, note that since $k_2$ is
generated by an ultrafilter ${\cal U}$ over $\gk$
and $2^\kappa=\kappa^+$ in both $V$ and $M$, we
know $|k_2(2^\gk)| = |k_2(\gk^+)| = |\{ f : f :
\gk \to \gk^+$ is a function$\}| =
|{[\gk^+]}^\gk| = \gk^+$. Thus, as $N[G_0]
\models ``|\wp(\FQ)| = k_2(2^\gk)$'', we can let
$\la D_\ga : \ga < \gk^+ \ra$ be an enumeration
in $V[G_0]$ of the dense open subsets of $\FQ$
present in $N[G_0]$. Since the $\gk$ closure of
$N$ with respect to either $M$ or $V$ implies
that the least element of the field of $\FQ$ is
above $\gk^+$, the definition of $\FQ$ as the
Easton support iteration which adds a
non-reflecting stationary set of ordinals of
cofinality $\eta$ to each $N$-measurable limit of
$N$-strong cardinals in the interval $(\gk,
k_2(\gk)]$ implies that $N[G_0] \models ``\FQ$ is
${\prec} \gk^+$-strategically closed''. Since the
standard arguments show that forcing with the
$\gk$-c.c$.$ partial ordering $\FP$ preserves
that $N[G_0]$ remains $\gk$-closed with respect
to either $M[G_0]$ or $V[G_0]$, we know that
$\FQ$ is ${\prec} \gk^+$-strategically closed in
both $M[G_0]$ and $V[G_0]$. We now construct
$G_1$ in either $M[G_0]$ or $V[G_0]$ as follows.
Fix a winning strategy for player II in the game
of length $\kappa^+$ for the partial ordering
$\FQ$ and use it to construct a play $\<q_\alpha
: \alpha<\kappa^+>$ of the game.
Since player II's moves are determined by her
strategy, we need only specify the moves of the
first player: if player II has just played the
condition $q_{2\alpha}$ at the (even) stage
$2\alpha$, let us direct player I to select and
then play a condition $q_{2\alpha+1}$ above
$q_{2\alpha}$ from the dense set $D_\alpha$.
Since the strategy plays at limit stages, this
completes the construction of the play
$\<q_\alpha : \alpha<\kappa^+>$. Let $G_1 = \{p
\in \FQ : \exists \ga < \gk^+\, (q_\ga \ge p)\}$
be the filter generated by this increasing
sequence of conditions. By construction, this
filter meets all the dense sets $D_\alpha$, and
so it is $N[G_0]$-generic over $\FQ$.

It remains to construct in $V[G_0]$ the desired
$N[G_0][G_1]$-generic object $G_2$ over $\FR$. To
do this, we first observe that as $M \models
``\gk$ is a measurable limit of strong
cardinals'', we can factor $k_1(\FP)$ as $\FP
\ast \dot \FS \ast \dot \FT$, where
$\forces_{\FP} ``\dot \FS = \dot \FP_{\eta,
\gk}$'',
and $\dot \FT$ is a term for the rest of
$k_1(\FP)$.

Note now that as in Lemma 2.4 of \cite{AC2}, $M
\models ``$No cardinal $\gd \in (\gk, \gl]$ is
strong''. Thus, the field of $\dot \FT$ is
composed of all $M$-measurable limits of
$M$-strong cardinals in the interval $(\gl,
k_1(\gk))$, which implies that in $M$,
$\forces_{\FP \ast \dot \FS} ``\dot \FT$ is
${\prec} \gl^+$-strategically closed''. Further,
since $\gl$ is a singular strong limit cardinal
of cofinality at least $\gk$, $|{[\gl]}^{< \gk}|
= \gl$. By Solovay's theorem \cite{So} that GCH
must hold at any singular strong limit cardinal
above a strongly compact cardinal, we know that
$2^\gl = \gl^+$. Therefore, as $k_1$ can be
assumed to be generated by an ultrafilter over
$P_\gk(\gl)$, we may calculate $|2^{k_1(\gk)}|^M
= |k_1(2^\gk)| = |k_1(\gk^+)| = |\{ f : f :
P_\gk(\gl) \to \gk^+$ is a function$\}| =
|{[\gk^+]}^\gl| = \gl^+$.

Work until otherwise specified in $M$.  Consider
the ``term forcing'' partial ordering $\FT^*$
(see \cite{F} for the first published account of
term forcing or \cite{C}, Section 1.2.5, page 8;
the notion is originally due to Laver) associated
with $\dot \FT$, i.e., $\tau \in \FT^*$
essentially iff $\tau$ is a term in the forcing
language with respect to $\FP \ast \dot \FS$ and
$\forces_{\FP \ast \dot \FS} ``\tau \in \dot
\FT$'', ordered by $\tau \ge \sigma$ iff
$\forces_{\FP \ast \dot \FS} ``\tau \ge
\sigma$''. Since this definition, taken
literally, would produce a proper class, we
restrict the terms appearing in it to a
sufficiently large set-sized collection (so that
any term $\tau$ forced by the trivial condition
to be in $\dot\FT$ will be forced by the trivial
condition to be equal to an element of $\FT^*$)
of size $k_1(\kappa)$ in $M$.\footnote{In the
official definition of $\FT^*$, the basic idea is
to include only the canonical terms. Since $\dot
\FT$ is forced to have cardinality $k_1(\kappa)$,
there is a set $\{\tau_\alpha :
\alpha<k_1(\kappa)\}$ of terms such that for any
other term $\tau$, if $\forces_{\FP \ast \dot
\FS} ``\tau \in \dot \FT$'', then there is a
dense set of conditions in $\FP \ast \dot \FS$
forcing $``\tau=\tau_\alpha$'' for various
$\alpha$. While this collection of terms may not
itself be adequate, we enlarge it as follows: for
each maximal antichain $A \subseteq \FP \ast \dot
\FS$ and each function $s: A \to \{\tau_\alpha :
\alpha<k_1(\kappa)\}$, there is (by the Mixing
Lemma of elementary forcing) a term $\tau_s$ such
that $p\forces``\tau_\ga=\tau_{s(p)}$'' for each
$p\in A$; let $\FT^*$ be the collection of all
such terms $\tau_s$, ranging over all maximal
antichains of $\FP \ast \dot \FS$. Since $\FP
\ast \dot \FS$ has size less than $k_1(\kappa)$
in $M$, the number of such terms is
$k_1(\kappa)$. And finally, if a term $\tau$ is
forced to be in $\dot \FT$, then elementary
forcing arguments establish that $\tau$ is forced
to be equal to $\tau_s$ for some $s$.} Since
$\forces_{\FP \ast \dot \FS} ``\dot \FT$ is
${\prec}\gl^+$-strategically closed'', it can
easily be verified that $\FT^*$ is also
${\prec}\gl^+$-strategically closed in $M$ and,
since $M^\gl \subseteq M$, in $V$ as well.

Since $M \models ``2^{k_1(\gk)} = {(k_1(\gk))}^+
= k_1(\gk^+)$'', we can let $\la D_\ga : \ga <
\gl^+ \ra$ be an enumeration in $V$ of the dense
open subsets of $\FT^*$ found in $M$ and argue as
we did when constructing $G_1$ to build in $V$ an
$M$-generic object $H_2$ over $\FT^*$. As readers
can verify for themselves, this line of reasoning
remains valid, in spite of the fact $\gl$ is
singular.

Note now that since $N$ is an ultrapower of $M$
via a normal ultrafilter ${\cal U} \in M$ over
$\gk$, Fact 2 of Section 1.2.2 of \cite{C} tells
us that $k_2 '' H_2$ generates an $N$-generic
object $G^*_2$ over $k_2(\FT^*)$. By
elementariness, $k_2(\FT^*)$ is the term forcing
in $N$ defined with respect to
$k_2(k_1(\FP_\gk)_{\gk + 1}) = \FP \ast \dot
\FQ$. Therefore, since $j(\FP) = k_2(k_1(\FP)) =
\FP \ast \dot \FQ \ast \dot \FR$, $G^*_2$ is
$N$-generic over $k_2(\FT^*)$, and $G_0 \ast G_1$
is $N$-generic over $k_2(\FP \ast \dot \FS)$, we
know by Fact 1 of Section 1.2.5 of \cite{C} (see
also \cite{F}) that $G_2=\{i_{G_0 \ast G_1}(\tau)
: \tau \in G^*_2\}$ is $N[G_0][G_1]$-generic over
$\FR$. Thus, in $V[G_0]$, the embedding $j : V
\to N$ lifts to $j : V[G_0] \to
N[G_0][G_1][G_2]$. This means that $V[G_0]
\models ``\gk$ is $\gl$ strongly compact''. As
$\gl$ was an arbitrary singular strong limit
cardinal of cofinality at least $\gk$, this
completes the proof of Lemma \ref{l1}.
\end{proof}

We remark that the proof of Lemma \ref{l1} shows
that a local version of this lemma is also
possible. Specifically, if $V \models ``{\rm GCH}
+ \gl \ge \gk$ is a cardinal + $\gk$ is a limit
of strong cardinals + $\gk$ is $\gl$
supercompact'', then $V^\FP \models ``\gk$ is
$\gl$ strongly compact''.

\begin{lemma}\label{l2}

$V^\FP \models ``\gk$ has trivial Mitchell
rank''.

\end{lemma}

\begin{proof}
Let $G$ be $V$-generic over $\FP$. If $V[G]
\models ``\gk$ does not have trivial Mitchell
rank'', then let $j : V[G] \to M[j(G)]$ be an
embedding generated by a normal measure over
$\gk$ in $V[G]$ witnessing this fact. In the
terminology of \cite{H1}, \cite{H2}, and
\cite{H3}, $\FP$ admits a gap below $\gk$, and so
by the Gap Forcing Theorem of \cite{H2} and
\cite{H3}, $j$ must lift an embedding $j : V \to
M$ that is definable in $V$. Since $j(\FP)$ also
admits a gap below $\kappa$ in $M$ and $\kappa$
is measurable in $M[j(G)]$, we similarly conclude
that $\kappa$ is measurable in $M$. Therefore,
since $\kappa$ is a limit of strong cardinals (we
have already noted that any supercompact cardinal
is a limit of strong cardinals), it follows that
$\gk$ is in the field of $j(\FP)$. Thus, there is
nontrivial forcing at stage $\kappa$, and so
$j(G) = G \ast S \ast H$, where $S$ is a
non-reflecting stationary set of ordinals added
by forcing over $M[G]$ with ${(\FP_{\eta,
\gk})}^{M[G]}$ at stage $\kappa$ and $H$ is
$M[G][S]$-generic for the rest of the forcing
$j(\FP)$. Since $V_{\kappa+1}\subseteq M\subseteq
V$, it follows that $V_{\gk + 1}^{V[G]} = V_{\gk
+ 1}^{M[G]}$. From this it follows that
${(\FP_{\eta, \gk})}^{M[G]} = {(\FP_{\eta,
\gk})}^{V[G]}$, and the dense open subsets of
what we can now unambiguously write as
$\FP_{\eta, \gk}$ are the same in both $M[G]$ and
$V[G]$. Thus, the set $S$, which is an element of
$V[G]$, is $V[G]$-generic over $\FP_{\eta, \gk}$,
a contradiction. This proves Lemma \ref{l2}.
\end{proof}

\begin{lemma}\label{l2a}

$V^\FP \models ``$No cardinal $\gg \in (\eta,
\gk)$ is strongly compact''.

\end{lemma}

\begin{proof}
By the definition of $\FP$, $V^\FP \models
``$Unboundedly many cardinals $\gg \in (\eta,
\gk)$ contain non-reflecting stationary sets of
ordinals of cofinality $\eta$''. Therefore, by
Theorem 4.8 of \cite{SRK} and the succeeding
remarks, $V^\FP \models ``$No cardinal $\gg \in
(\eta, \gk)$ is strongly compact''. This proves
Lemma \ref{l2a}.
\end{proof}

The proofs of Lemmas \ref{l1} - \ref{l2a}
complete the proof of Proposition \ref{p1}.
\end{proof}

We turn now to the proof of Proposition \ref{p2}.

\begin{proof}
Let $V \models ``$ZFC + $\gk$ is supercompact'',
with $\eta < \gk$ fixed but arbitrary. As in the
proof of Proposition \ref{p1}, we can assume
without loss of generality that $V \models
``2^\gk = \gk^+$'' and that if necessary, this
has been forced by the use of an $\eta$-directed
closed partial ordering having cardinality
$2^\gk$ that preserves all supercompact cardinals
above $\eta$ and creates no new supercompact or
measurable cardinals above $\eta$. The proof of
Proposition \ref{p2} is then given by the
following two lemmas.

\begin{lemma}\label{n1}

There is an $\eta$-directed closed partial
ordering $\FP$, preserving all supercompact
cardinals above $\eta$ except for $\gk$, so that
$V^\FP \models ``\gk$ is a non-supercompact
strongly compact cardinal which contains a club
disjoint from the measurable cardinals''.
Furthermore, every strongly compact cardinal in
$V^{\FP}$ in the interval $(\eta, \gk]$ is either
supercompact in $V$ or a measurable limit of
supercompact cardinals in $V$.

\end{lemma}

\begin{proof}
This proof follows the main idea of
\cite[Corollary 4.3]{H4}, suitably modified so as
to ensure the requirement concerning strongly
compact cardinals in the extension. Let $f$ be a
uniform Laver function for all supercompact
cardinals in the interval $(\eta, \gk]$ (the
existence of such a function is shown in
\cite{KM} and \cite{A98}). The forcing $\FP^0$
will be a kind of modified lottery preparation
with respect to $f$. Specifically, $\FP^0$ is the
Easton support iteration of length $\kappa$ which
begins by adding a Cohen subset to $\eta^+$.
$\FP^0$ then has nontrivial forcing at stage
$\ga$ only when $\ga > \eta$ is a measurable
cardinal, $\ga \in \dom(f)$, $f(\ga)$ is an
ordinal, $\ga \le f(\ga)$
and $f''\ga \subseteq V_\ga$. At such stages, the
forcing is the lottery sum in $V^{\FP_\ga}$ of
all partial orderings in $H(f(\alpha)^+)$ having
$\beta$-directed closed dense subsets for every
$\beta<\alpha$. After the lottery sum forcing, we
perform the forcing to add a non-reflecting
stationary subset to the next inaccessible
cardinal above $f(\alpha)$, consisting of
ordinals of cofinality $\max(\eta^+,
\beta_\alpha^+)$, where $\beta_\alpha$ is the
supremum of the supercompact cardinals of $V$
below $\alpha$. (Note: this will prevent any
cardinals in the interval $(\max(\eta,
\beta_\alpha) ,f(\alpha)]$ from becoming strongly
compact in the extension.) We may assume that
$\dom(f)$ contains no supercompact cardinals, so
that the forcing at any supercompact cardinal
stage is trivial.

Let us argue that forcing with $\FP^0$ preserves
all supercompact cardinals $\delta$ of $V$ in the
interval $(\eta,\kappa]$. For this, it will
suffice for us to argue that $\delta$ is
supercompact in $V^{\FP_\delta}$ and
indestructible there by any further
$\delta$-directed closed forcing. This suffices
because there is no forcing right at stage
$\delta$ (as it is not in the domain of $f$), and
the subsequent forcing from stage $\delta$ up to
$\kappa$ has a $\delta$-directed closed dense
subset. Following \cite[Corollary 4.6]{H4}, fix
any $\delta$-directed closed forcing $\FQ \in
V^{\FP_\delta}$ and any $\lambda\geq\delta$.
Choose $\theta > {(\max(2^{\lambda^{<\delta}},
|{\rm TC}(\dot\FQ)|))}^+$. Since $f$ is a Laver
function, there is a $\theta$ supercompactness
embedding $j:V\to M$ with critical point $\delta$
and $j(f)(\delta)=\theta$. Let $G_\delta*g
\subseteq \FP_\delta*\dot\FQ$ be $V$-generic.
Since $\FQ$ is allowed in the stage $\delta$
lottery of $j(\FP_\delta)$, we may work above the
condition opting for this partial ordering, and
factor the forcing as $j(\FP_\delta)=\FP_\delta
\ast \dot \FQ \ast \dot \FP_\delta^{j(\delta)}$,
where $\dot \FP_\delta^{j(\delta)}$ is a term for
the rest of the forcing from stages $\delta$ up
to $j(\delta)$, starting with the non-reflecting
stationary set forcing at stage $\delta$. Force
to add $G_\delta^{j(\delta)} \subseteq
\FP_\delta^{j(\delta)}$, and lift the embedding
to $j:V[G_\delta]\to M[j(G_\delta)]$ in
$V[G_\delta][g][G_\delta^{j(\delta)}]$, where
$j(G_\delta)=G_\delta*g*G_\delta^{j(\delta)}$.
Using a master condition above $j''g$, similarly
force to add $j(g) \subseteq j(\FQ)$, and lift
the embedding to $j:V[G_\delta][g]\to
M[j(G_\delta)][j(g)]$ in
$V[G_\delta][g][G_\delta^{j(\delta)}][j(g)]$. The
point is now that the induced $\lambda$
supercompactness measure $\mu$, defined by $X\in
\mu$ iff $j''\lambda\in j(X)$, has size
$2^{\lambda^{{<}\delta}}$, and therefore $\mu$
could not have been added by the extra forcing
$\FP_\delta^{j(\delta)}*j(\dot\FQ)$, since that
forcing is $2^{{\gl}^{< \gd}}$-strategically
closed. Hence, the measure $\mu$ must be in
$V[G_\delta][g]$, and so $\delta$ is $\lambda$
supercompact there, as desired. So every
supercompact cardinal of $V$ in the interval
$(\eta,\kappa]$ is preserved to $V^{\FP^0}$ and
becomes indestructible there.

Now let $\FC_\gk$ be the partial ordering defined
in $V^{\FP^0}$ which adds a club of
non-measurable cardinals to $\gk$ above $\eta$,
i.e., $\FC_\gk = \{c : c$ is a closed, bounded
subset of $\gk$ containing no cardinals that are
measurable in $V^{\FP^0}$ and $\eta<\sup(c)\}$,
ordered by end-extension. Let $\FP = \FP^0 \ast
\dot \FC_\gk$. For every $\gd < \gk$, the set of
elements of $\FC_\gk$ which mention an element
above $\gd$ is a $\gd^+$-directed closed dense
open subset of $\FC_\gk$. This can be seen by
taking the union of a $\gd$-chain of these sorts
of conditions and adding the supremum (which
cannot be measurable because it is not regular).
Thus, the measurable cardinals below $\gk$ in
$V^{\FP^0}$ and $V^{\FP^0 \ast \dot \FC_\gk} =
V^\FP$ are the same. Also, $\FP$ is
$\eta$-directed closed, and for any
$V$-supercompact cardinal $\gd \in (\eta, \gk]$,
by indestructibility, $V^{\FP^0 \ast \dot
\FC_\gk} = V^\FP \models ``\gd$ is a supercompact
cardinal whose supercompactness is indestructible
under $\gd$-directed closed forcing''. In
addition, since by its definition, $|\FP| = \gk$,
the results of \cite{LS} imply that forcing with
$\FP$ preserves all supercompact cardinals above
$\gk$, and $V^\FP \models ``2^\gk = \gk^+$''.

Since $f$ is a Laver function for $\gk$, we know
that for any cardinal $\gl \ge \gk$, there is a
supercompact ultrafilter ${\cal U}_0$ over
$P_\gk(\gl)$ so that for $j_0 : V \to M$ the
associated elementary embedding generated by
${\cal U}_0$, $j_0(f)(\gk) = \gl$. As $M \models
``{({[{\rm id}]}_{{\cal U}_0})}^{M} = \la
j_0(\ga) : \ga < \gl \ra$'', it follows that $M
\models ``{|{[{\rm id}]}_{{\cal U}_0}|}^{M} =
\gl$''. Therefore, $M \models ``j_0(f)(\gk) \ge
{|{[{\rm id}]}_{{\cal U}_0}|}^{M}$''. Further, if
$\gl \ge 2^\gk$, as in the proof of Lemma
\ref{l1}, we can find $j_1 : M \to N$ an
elementary embedding generated by a normal
ultrafilter ${\cal U}_1 \in M$ of trivial
Mitchell rank so that $N \models ``\gk$ isn't
measurable''. Let $\cal U$ be the $\gk$-additive,
fine measure over $P_\gk(\gl)$ defined by $x \in
{\cal U}$ iff $j_1(\la j_0(\ga) : \ga < \gl \ra)
\in x$, with the associated elementary embedding
$j : V \to M^*$. Thus, $\gk \in j(\{\gd < \gk :
\gd$ isn't measurable$\})$ and $j(f)(\gk) \ge
{|{[{\rm id}]}_{{\cal U}}|}^{M^*}$.  We are
therefore in a position to apply the argument of
\cite[Theorem 4.2]{H4} to conclude that $V^{\FP^0
\ast \dot \FC_\gk} = V^\FP \models ``\gk$ is
strongly compact''. And we've explicitly added a
club of non-measurable cardinals to $\kappa$.

It remains to check that every strongly compact
cardinal of $V^\FP$ in the interval $(\eta, \gk]$
is either supercompact in $V$ or a measurable
limit of supercompact cardinals in $V$. Note that
by construction, whenever $\gamma \in (\eta,
\gk]$ is a supercompact cardinal of $V$ which
isn't a limit of supercompact cardinals, then for
unboundedly many $\lambda$ between $\max(\eta,
\gb_\gg)$ and $\gg$,
we have added a non-reflecting stationary subset
to $\lambda$ consisting of ordinals of cofinality
$\max(\eta^+, \gb^+_\gg)$. This necessarily
destroys all strongly compact cardinals between
$\max(\eta, \gb_\gg)$ and $\gg$. So no strongly
compact cardinals in the extension in the
interval $(\eta, \gk]$ lie between the
supercompact cardinals of $V$ in the interval
$(\eta, \gk]$. Since furthermore by the Gap
Forcing Theorem of \cite{H2} and \cite{H3}, no
new measurable or supercompact cardinals were
created, we conclude that every strongly compact
cardinal in the extension is either supercompact
in $V$ or a measurable limit of supercompact
cardinals in $V$, as we claimed. This completes
the proof of Lemma \ref{n1}.
\end{proof}

\begin{lemma}\label{n2}

There is an $\eta$-directed closed partial
ordering $\FP$, preserving all supercompact
cardinals above $\eta$ except for $\gk$, so that
$V^\FP \models ``\gk$ is a non-supercompact
strongly compact strong cardinal''. Furthermore,
every strongly compact cardinal in $V^\FP$ in the
interval $(\eta, \gk]$ is either supercompact in
$V$ or a measurable limit of supercompact
cardinals in $V$.

\end{lemma}

\begin{proof}
Let $\FP$ be the Easton support iteration of
length $\gk$ which begins by adding a Cohen
subset to $\eta^+$. $\FP$ then has nontrivial
forcing only at those stages $\ga > \eta$ which
are strong cardinal limits of strong cardinals.
At such a stage $\ga$, we force with the lottery
sum of all $\ga$-directed closed partial
orderings having rank below the least strong
cardinal $\gd$ above $\ga$, which add a Cohen
subset to $\ga$. We next perform the forcing to
add a non-reflecting stationary subset of
ordinals of cofinality $\max(\eta^+, \gb^+_\ga)$
to $\gd$, where $\gb_\ga$ is as in Lemma
\ref{n1}. It follows easily that $\FP$ is
$\eta$-directed closed.

For any $\gl > \gk$ so that $\gl = \beth_\gl$, we
can choose $j : V \to M$ to be an elementary
embedding witnessing the $\gl$ strongness of
$\gk$ so that $M \models ``\gk$ isn't $\gl$
strong''. This means that by the definition of
$\FP$, no forcing is done in $M$ at stage $\gk$.
Therefore, the standard lifting argument for
strongness embeddings will show that $\kappa$
remains strong in $V^\P$.\footnote{One factors
through by the normal measure, constructs a
generic over the normal measure ultrapower, and
pushes this generic up to the strongness
ultrapower. See, for example, Lemma 2.5 of
\cite{AC2}. A key part of the argument is that as
in Lemma 2.5 of \cite{AC2}, since no cardinal
$\gd \in [\gk, \gl]$ in the strongness ultrapower
is strong, the first stage of nontrivial forcing
in the strongness ultrapower takes place well
after stage $\gl$.} Further, if now $\gl \ge
2^\gk$ and $j : V \to M$ is an elementary
embedding witnessing the $\gl$ supercompactness
of $\gk$, then by remarks made in the proof of
Lemma \ref{l1}, $\gk$ is a strong cardinal limit
of strong cardinals in both $V$ and $M$. This
means that in $M$, the forcing for adding a Cohen
subset to $\gk$ is part of the lottery sum found
at stage $\gk$ in the definition of $j(\FP)$.
Hence, since we are able to opt for this forcing
at stage $\gk$ in $M$, we can apply the argument
given in the proof of Lemma \ref{l1} to show that
$V^\FP \models ``\gk$ is strongly compact''. In
addition, we can modify the proof of Lemma
\ref{l2} by replacing the embedding $j$ with an
embedding that is alleged to witness the fact
that $\gk$ is $2^\gk$ supercompact
in the generic extension and by replacing the
non-reflecting stationary set of ordinals of
Lemma \ref{l2} with a Cohen subset of $\gk$. This
last change is possible since the stage $\kappa$
forcing in $M$ must add a Cohen subset to
$\kappa$. The proof of Lemma \ref{l2} now goes
through in an analogous manner as earlier to show
that after forcing with $\FP$, $\gk$ isn't
$2^\gk$ supercompact.
Therefore, since $|\FP| = \gk$, we know $V^\FP
\models ``2^\gk = \gk^+$'', and the results of
\cite{LS} once again show that forcing with $\FP$
preserves all supercompact cardinals above $\gk$.

Since the proof that every strongly compact
cardinal in $V^\FP$ in the interval $(\eta, \gk]$
is either supercompact in $V$ or a measurable
limit of supercompact cardinals in $V$ is exactly
as given in the proof of Lemma \ref{n1}, we
complete the proof of Lemma \ref{n2} by showing
that $\FP$ preserves all supercompact cardinals
in the interval $(\eta, \gk)$. To see this, let
$\gd \in (\eta, \gk)$ be supercompact. Let $\gl >
\gk$ be a singular strong limit cardinal of
cofinality $\gd$ (so the GCH holds at $\lambda$).
Choose $j : V \to M$ an elementary embedding
witnessing the $\gl$ supercompactness of $\gd$ so
that $M \models ``\gd$ isn't $\gl$
supercompact''. Since $\gd$ is a strong cardinal
limit of strong cardinals, $\gd$ is a stage in
the definition of $\FP$ at which a nontrivial
forcing is done, i.e., if we write $\FP = \FP_\gd
\ast \dot \FP^\gd$, $\dot \FP^\gd$ will be a term
for a partial ordering that is $\gd$-directed
closed and adds a Cohen subset to $\gd$. Because
$\delta$ is ${<}\lambda$ supercompact but not
$\lambda$ supercompact in $M$, it follows as
before that $M \models ``$There are no strong
cardinals in the interval $(\gd, \gl]$''. This
means that $\FP^\gd$ is an allowable choice in
the stage $\gd$ lottery in $M^{\FP_\gd}$, and any
further nontrivial forcing in $M$ takes place
well after stage $\gl$. Therefore, if $G_0$ is
$V$-generic over $\FP_\gd$ and $G_1$ is
$V[G_0]$-generic over $\FP^\gd$, in
$V[G_0][G_1]$, standard arguments show that $j$
lifts to $j : V[G_0][G_1] \to
M[G_0][G_1][G_2][G_3]$, where $G_2$ and $G_3$ are
suitably generic objects constructed in
$V[G_0][G_1]$, and $G_3$ contains a master
condition for $G_1$. We can thus infer that
$V^\FP \models ``\gd$ is $\gl$ supercompact''.
Since $\gl$ can be chosen to be arbitrarily
large, this completes the proof of Lemma
\ref{n2}.
\end{proof}

Lemmas \ref{n1} and \ref{n2} complete the proof
of Proposition \ref{p2}.
\end{proof}

We conclude Section \ref{s2} by remarking that it
is possible to modify the definition of $\FP$
given in the proofs of Propositions \ref{p1} and
\ref{p2} so that after forcing with $\FP$, $\gk$
retains a nontrivial degree of supercompactness.
To do this, $\FP$ is first altered so as
initially to force $2^\gk = \gk^+$ and $2^{\gk^+}
= \gk^{++}$ if necessary via an $\eta$-directed
closed partial ordering that preserves all
$V$-supercompact cardinals, while admitting a gap
below the least inaccessible above $\eta$. The
forcing $\FP$ is then defined as in the proof of
Proposition \ref{p1} and Lemma \ref{n2}, except
that nontrivial forcing is done only at stages
$\ga$ which are $\ga^+$ supercompact and are a
limit of strong cardinals in the partial ordering
which is the analogue of the one defined in
Proposition \ref{p1}, or at stages $\ga$ which
are $\ga^+$ supercompact and are strong cardinal
limits of strong cardinals in the partial
ordering which is the analogue of the one defined
in Lemma \ref{n2}. If $j$ in the proof of Lemma
\ref{l2} is chosen as a $\gk^{+}$
supercompactness embedding witnessing that $\gk$
has nontrivial Mitchell rank with respect to
$\gk^+$ supercompactness (meaning that there is a
$\kappa^+$ supercompactness embedding $j : V \to
M$ with $M \models ``\gk$ is $\gk^+$
supercompact'') instead of an embedding
witnessing that $\gk$ has nontrivial Mitchell
rank and the word ``measurable'' is replaced by
the phrase ``$\gk^+$ supercompact'', then the
remainder of the proof of Lemma \ref{l2} suitably
modified shows that $V^\FP \models ``\gk$ has
trivial Mitchell rank with respect to $\gk^+$
supercompactness''. By replacing the embedding
$k_2$ in Lemma \ref{l1} with a $\gk^+$
supercompactness embedding so that $M \models
``\gk$ isn't $\gk^+$ supercompact'', the proof of
Lemma \ref{l1} also goes through with slight
modifications and shows that $V^\FP \models
``\gk$ is both strongly compact and $\gk^+$
supercompact''. The proof that all
$V$-supercompact cardinals above $\eta$ except
for $\gk$ remain supercompact is virtually
identical to the one given in Lemma \ref{n2} in
the appropriate analogue of Lemma \ref{n2}, and
the proofs that there are no strongly compact
cardinals in the interval $(\eta, \gk)$ in the
appropriate analogue of Lemma \ref{l2a} or that
the strongly compact cardinals in $V^\FP$ in the
interval $(\eta, \gk]$ are either supercompact in
$V$ or measurable limits of supercompact
cardinals in $V$ in the appropriate analogue of
Lemma \ref{n2} are identical to the ones given in
the proofs of these lemmas.
Thus, by modifying the definition of $\FP$, it is
possible to produce an $\eta$-directed closed
partial ordering $\FP$ with the same properties
as in Propositions \ref{p1} and \ref{p2} except
that in $V^\FP$, $\gk$ is a non-supercompact
strongly compact $\kappa^+$ supercompact cardinal
having trivial Mitchell rank with respect to
$\gk^+$ supercompactness (see also pages 113--114
of \cite{A97}). Further modifications allow
$V^\FP$ to witness even larger degrees of
supercompactness for $\gk$, while remaining
strongly compact and not supercompact. One can
also arrange that in $V^{\FP}$ the cardinal
$\kappa$ has trivial Mitchell rank with respect
to its degree of supercompactness.

\section{The Proof of the Main Theorem}\label{s3}

We turn now to the proof of our Main Theorem,
Theorem \ref{t1}.

\begin{proof}
Let $V_0 \models ``$ZFC + $\K$ is the class of
supercompact cardinals + ${\cal A} \subseteq \K$
and $\cal A$ contains none of its limit points.''
By initially forcing with an Easton support
iteration $\FP^*$ that first adds a Cohen subset
to $\go$, next forces GCH if necessary at all
measurable cardinals $\gk$ by adding a Cohen
subset to $\gk^+$ (which preserves all
supercompact cardinals and creates no new
supercompact or measurable cardinals) and then is
followed by a modified version of the partial
ordering used in the proof of Theorem 1 of
\cite{A98}, we may assume that $V = V^{\FP^*}_0
\models ``$ZFC + $\K$ is the class of
supercompact cardinals + $2^\gk = \gk^+$ if $\gk$
is supercompact + Every supercompact cardinal
$\gk$ is Laver indestructible \cite{L} under
$\gk^+$-directed closed forcing + The strongly
compact and supercompact cardinals coincide
precisely, except possibly at measurable limit
points''. The modification is to allow only
$\gd^+$-directed closed partial orderings in the
indestructibility forcing at stage $\gd$. In
addition to forcing the previous properties, this
modification will ensure that all measurable
limits of supercompact cardinals at which GCH
holds in $V_0$ are preserved and continue to
satisfy GCH in $V$.\footnote{Let us outline the
proof that $\FP^*$ accomplishes this. The portion
of $\FP^*$ that forces GCH at a $V_0$-measurable
cardinal $\gk$ is an iteration that can be
factored as $\FQ_0 \ast \dot \FQ_1$, where the
field of $\FQ_0$ is composed of ordinals below
$\gk$. If $\gk$ is any measurable cardinal at
which GCH holds, then by choosing $j : V \to M$
as an elementary embedding witnessing $\gk$'s
measurability so that $M \models ``\gk$ isn't
measurable'', we can show via a standard lifting
argument (such as given in the proof of Theorem
3.5 of \cite{H4}) that $j$ lifts to an elementary
embedding witnessing $\gk$'s measurability after
forcing with $\FQ_0$. Since $\forces_{\FQ_0}
``\dot \FQ_1$ is $\gk^+$-directed closed'', $\gk$
remains measurable after forcing with $\FQ_1$.
Further, since $|\FQ_0| \le \gk$, GCH at $\gk$ is
preserved after forcing with $\FQ_0 \ast \dot
\FQ_1$. Then, since the remainder of $\FP^*$ is a
Laver style iteration for forcing
indestructibility which at each nontrivial stage
$\gd$ does a forcing which is $\gd^+$-directed
closed, it can be factored as $\FQ_2 \ast \dot
\FQ_3$, where for $\gk$ a $V_0$-measurable limit
of supercompact cardinals, the field of $\FQ_2$
is composed of ordinals below $\gk$, $|\FQ_2| \le
\gk$ and $\forces_{\FQ_2} ``\dot \FQ_3$ is
$\gk^+$-directed closed''. The argument that
forcing with $\FQ_0 \ast \dot \FQ_1$ preserves
both $\gk$'s measurability and GCH at $\gk$ can
now be applied to $\FQ_2 \ast \dot \FQ_3$.
Standard arguments for the preservation of
supercompactness under Easton support iterations
show that any supercompact cardinal is preserved
after forcing with $\FQ_0 \ast \dot \FQ_1$, and
the arguments of \cite{A98} show that the
remaining claimed properties of $\FP^* = \FQ_0
\ast \dot \FQ_1 \ast \dot \FQ_2 \ast \dot \FQ_3$
hold. The Gap Forcing Theorem of \cite{H2} and
\cite{H3} then shows that $\K$ is precisely the
class of supercompact cardinals in $V$.}

We now define the forcing $\FP$ that will
accomplish our goals. Work in $V$. For each
$\kappa\in\cal A$, let
$\delta_\kappa = 2^{\sup (\cal
A\intersect\kappa)}$, with $\delta_{\kappa_0} =
\ha_2$ for $\gk_0$ the least element of ${\cal
A}$. Since ${\cal A}$ contains none of its limit
points, it follows that $\gd_\gk < \gk$. For each
$\kappa\in\cal A$, let $\FP_\gk$ be the partial
ordering which
forces $\gk$ to be a non-supercompact strongly
compact cardinal by first taking $\eta = \gd_\gk$
and then using any of the $\eta =
\gd_\gk$-directed closed partial orderings given
in Proposition \ref{p2}. Let $\FP$ be the Easton
support product $\prod_{\gk \in {\cal A}}
\FP_\gk$. Please note that this is a product, not
an iteration. The field of $\FP_\kappa$ lies in
the interval $(\delta_\kappa,\kappa)$, which
contains no elements of $\cal A$, and so the
fields of the partial orderings $\FP_\kappa$
occur in disjoint blocks. Although this may be
class forcing, the standard Easton arguments show
$V^{\FP} \models {\rm ZFC}$.

\begin{lemma}\label{l3}

If $\gk \in \K - {\cal A}$, then $V^{\FP} \models
``\gk$ is supercompact''.

\end{lemma}

\begin{proof}
Suppose that $\gk \in \K - {\cal A}$. Let $\eta$
be the least element of ${\cal A}$ above $\gk$,
and factor the forcing in $V$ as $\FP = \FQ^\eta
\times \FP_\eta \times \FQ_{< \eta}$, where
$\FQ^\eta = \prod_{\gb > \eta, \gb \in {\cal A}}
\FP_\gb$ and $\FQ_{< \eta} = \prod_{\gb < \eta,
\gb \in {\cal A}} \FP_\gb$. Please observe that
$\FQ^\eta$ is $\eta^+$-directed closed.

Suppose that $\kappa$ is not a limit point of
$\cal A$. It follows that $\cal A$ is bounded
below $\kappa$ and so $\delta_\eta<\kappa$ and
$|\FQ_{<\eta}|<\kappa$. Since $\kappa$'s
supercompactness is indestructible in $V$ and
$\FQ^\eta$ is $\eta^+$-directed closed, we know
that $\kappa$ is supercompact in $V^{\FQ^\eta}$.

We claim that the strong cardinals below $\eta$
are the same in $V^{\FQ^\eta}$ as in $V$. To see
this, observe first that $\eta$ is supercompact
and therefore strong in both $V$ and, by
indestructibility, in $V^{\FQ^\eta}$. Therefore,
a cardinal $\delta<\eta$ is strong in either $V$
or $V^{\FQ^\eta}$ iff it is $\sigma$ strong for
every $\sigma<\eta$, since by the second
paragraph of Lemma 2.1 of \cite{AC2}, if $\ga$ is
$\gb$ strong for every $\gb < \gg$ and and $\gg$
is strong, then $\ga$ is strong. But since
$\FQ^\eta$ is $\eta^+$-directed closed, it
neither creates nor destroys any extenders below
$\eta$, and so the two models agree on strongness
below $\eta$.

In addition, since the models agree up to $\eta$,
the construction of a universal Laver function
for the supercompact cardinals in the interval
$(\delta_\eta,\eta]$ is the same in $V$ or
$V^{\FQ^\eta}$. Therefore, the forcing $\FP_\eta$
satisfies the same definition in either $V$ or
$V^{\FQ^\eta}$, and so by applying Proposition
\ref{p2} in $V^{\FQ^\eta}$, since $\FP_\eta$
preserves all supercompact cardinals in the
interval $(\delta_\eta,\eta)$, we conclude
$V^{\FQ^\eta \times \FP_\eta} \models ``\gk$ is
supercompact''. Finally, since $|\FQ_{{<} \eta}|
< \gk$, we conclude by \cite{LS} that
$V^{\FP}=V^{\FQ^\eta \times \FP_\eta
\times\FQ_{{<} \eta}} \models``\kappa$ is
supercompact'', as desired.

Assume next that $\kappa\in\cal K-A$ and $\kappa$
is a limit point of $\cal A$. In this case,
$\gd_\eta = 2^\gk$, and so the partial ordering
$\FQ^\eta\times\FP_\eta$ is
$\delta_\eta$-directed closed. By
indestructibility, therefore, $\kappa$ is
supercompact in the model $\ov
V=V^{\FQ^\eta\times\FP_\eta}$. Furthermore, $V
\models ``|\FQ_{< \eta}| = \gk$''.

Choose any cardinal $\gl > \gk$ and let $\gg =
|2^{{\gl}^{< \gk}}|$. Take $j : \ov V \to M$ to
be an elementary embedding witnessing the $\gg$
supercompactness of $\gk$ in $\ov V$ so that $M
\models ``\gk$ isn't $\gg$ supercompact''. It
must then be the case, as in Lemma \ref{l1} and
Lemma 2.4 of \cite{AC2}, that $M \models ``$No
cardinal $\gd \in (\gk, \gg]$ is strong''.
Writing $j(\FQ_{< \eta})$ as $\FQ_{< \eta} \times
\FQ^*$, this means that the least ordinal in the
field of $\FQ^*$ is above $\gg$. Thus, if $G$ is
$\ov V$-generic over $\FQ_{< \eta}$ and $H$ is
$\ov V[G]$-generic over $\FQ^*$, in $\ov V[G
\times H]$, $j$ lifts to $\ov j : \ov V[G] \to
M[G \times H]$ via the definition $\ov
j(i_G(\tau)) = i_{G \times H}(j(\tau))$. Since $M
\models ``\FQ^*$ is $\gg$-strategically closed''
and $M^\gg \subseteq M$, it follows that for any
cardinal $\gd \le \gg$, the two models $\ov V[G]$
and $\ov V[G \times H] = \ov V[H \times G]$
contain the same subsets of $\gd$. This means the
supercompactness measure ${\cal U}$ over
${(P_{\gk}(\gl))}^{\ov V[G]}$ in $\ov V[G \times
H]$ given by $x \in {\cal U}$ iff $\la j(\gb) :
\gb < \gl \ra \in \ov j(x)$ is in $\ov V[G]$.
Hence, $V^{\FP}=V^{\FQ^\eta \times \FP_\eta
\times \FQ_{< \eta}}
 \models ``\gk$ is supercompact'', as desired.
This completes the proof of Lemma \ref{l3}.
\end{proof}

\begin{lemma}\label{jdh0}

If $\gk \in {\cal A}$, then $V^{\FP} \models
``\gk$ is strongly compact but not supercompact +
$2^\gk = \gk^+$''.

\end{lemma}

\begin{proof}
In analogy with the above argument, we factor the
forcing as $\FP = \FQ^\gk \times \FP_\gk \times
\FQ_{< \gk}$, where $\FQ^\gk = \prod_{\gb
> \gk, \gb \in {\cal A}} \FP_\gb$, and $\FQ_{< \gk} = \prod_{\gb < \gk,
\gb \in {\cal A}} \FP_\gb$. Since $\FQ^\gk$ is
$\gk^+$-directed closed, we know once again by
indestructibility that $\kappa$ is supercompact
in $V^{\FQ^\kappa}$. We also know again that the
strong cardinals below $\kappa$ of $V$ and
$V^{\FQ^\kappa}$ are the same and that the
partial ordering $\FP_\kappa$ of Proposition
\ref{p2} is constructed in the same way in $V$ as
in $V^{\FQ^\kappa}$. Therefore, by Proposition
\ref{p2}, we know that $\kappa$ is a
non-supercompact strongly compact cardinal in
$V^{\FQ^\kappa\times\FP_\kappa}$ (and either has
a club disjoint from the measurable cardinals or
else remains a strong cardinal, depending on the
version of $\FP_\kappa$ selected). Finally, since
our assumption on ${\cal A}$ guarantees that
$\kappa$ is not a limit point of $\cal A$, we
know $|\FQ_{{<}\kappa}| < \gk$. By the results of
\cite{LS} and \cite{HW}, therefore, $V^{\FQ^\gk
\times \FP_\gk \times \FQ_{{<} \gk}} = V^\FP
\models ``\gk$ is a non-supercompact strongly
compact cardinal which either has a club disjoint
from the measurables or is a strong cardinal''.
Lastly, we observe that none of the three factors
destroys $2^\kappa=\kappa^+$, so the proof of
Lemma \ref{jdh0} is complete.
\end{proof}

\begin{lemma}\label{jdh1}

If $V \models ``\gk$ is a measurable limit of
${\cal A}$ + $2^\gk = \gk^+$'', then $V^\FP
\models ``\gk$ is a measurable limit of ${\cal
A}$ + $2^\gk = \gk^+$''.

\end{lemma}

\begin{proof}
Suppose that $V \models ``\gk$ is a measurable
limit of ${\cal A}$ + $2^\gk = \gk^+$''.
By hypothesis, we know that $\kappa\notin\cal A$.
We may therefore factor the forcing $\FP$ as $\FP
= \FQ^\gk \times \FQ_{< \gk}$, where $\FQ^\gk =
\prod_{\gb > \gk, \gb \in {\cal A}} \FP_\gb$, and
$\FQ_{< \gk} = \prod_{\gb < \gk, \gb \in {\cal
A}} \FP_\gb$. Fix an elementary embedding $j:V\to
M$ arising from a normal measure over $\kappa$,
such that $M \models ``\gk$ isn't measurable''.

First, we observe that the forcing $\FQ^\kappa$
is $\kappa^+$-directed closed, and therefore
cannot destroy the measurability of $\kappa$. So
we need only prove that $\kappa$ is measurable
after forcing with $\FQ_{<\kappa}$.
Because $\kappa$ is a limit point of $\cal A$,
the forcing $j(\FQ_{<\kappa})$ factors as
$\FQ_{<\kappa}\times\FQ_{\kappa,j(\kappa)}$,
where $\FQ_{\kappa,j(\kappa)}$ is the product of
the partial orderings $\FP_\beta$ in $M$ for
$\beta\in j({\cal A}) \cap [\kappa,j(\kappa))$.
Since $\kappa$ is not measurable in $M$, it is
definitely not in $j(\cal A)$, and so the first
element of the field of $\FQ_{\kappa,j(\kappa)}$
is strictly above $\kappa$. It follows that this
forcing is $\kappa^+$-directed closed in $M$, and
so since $2^\kappa=\kappa^+$ in $V$, we may
employ the usual diagonalization argument (see,
e.g., the construction of the generic object
$G_1$ in Lemma \ref{l1}) to build in $V$ an
$M$-generic filter $G_{\kappa,j(\kappa)}\subseteq
\FQ_{\kappa,j(\kappa)}$. Putting this together
with any $V$-generic $G_{<\kappa}$ for
$\FQ_{<\kappa}$, we may lift the embedding to
$j:V[G_{<\kappa}]\to M[j(G_{<\kappa})]$, where
$j(G_{<\kappa})=G_{<\kappa}\times
G_{\kappa,j(\kappa)}$. This lifted embedding
witnesses the measurability of $\kappa$ in
$V[G_{<\kappa}]$, as desired. Finally, since
$\FQ^\gk$ is $\gk^+$-directed closed and
$|\FQ_{{<} \gk}| \le \gk$, $V^{\FQ^\gk \times
\FQ_{{<} \gk}} = V^\FP \models ``2^\gk =
\gk^+$''. This proves Lemma \ref{jdh1}.
\end{proof}

\begin{lemma}\label{l4}

The strongly compact cardinals of $V^\FP$ are
precisely the cardinals of $\cal K$ and their
measurable limit points, and these are all
strongly compact in $V$ and $V_0$. In addition,
the supercompact cardinals of $V^\FP$ are all
supercompact in $V$ and $V_0$ as well.

\end{lemma}

\begin{proof}
We have already proved that the cardinals of
$\cal K$ remain strongly compact in $V^\FP$, so
suppose towards a contradiction that $V^{\FP}
\models ``\gd \not\in \K$ is strongly compact and
isn't a measurable limit point of $\K$''. If
$\cal A$ is bounded below $\delta$, then $|\FP| <
\gd$, and so by the results of \cite{LS},
$\delta$ is strongly compact in $V$, contrary to
our assumption that the strongly compact
cardinals of $V$ are either in $\cal K$ or
measurable limits of $\cal K$. So we may assume
$\cal A$ is not bounded below $\delta$. Since
$\delta$ cannot be a limit point of $\cal A$,
there is a least element $\gk$ in ${\cal A}$
above $\gd$, and $\gd \in (\gd_\gk, \gk)$.

As in Lemma \ref{l3}, factor $\FP$ as $\FQ^\gk
\times \FP_\gk \times \FQ_{< \gk}$. Once again,
regardless of which version is chosen, $\FP_\gk$
is constructed the same in either $V$ or
$V^{\FQ^\kappa}$, and has the properties
identified in Proposition \ref{p2} in either
model. We therefore know that in particular,
since $\gd$ isn't in $V$ either an element of
$\K$ or a measurable limit of elements of $\K$,
$V^{\FQ^\gk \times \FP_\gk} \models ``\gd$ isn't
strongly compact''. Hence, once again, the
results of \cite{LS} yield that $V^{\FQ^\gk
\times \FP_\gk \times \FQ_{< \gk}} = V^{\FP}
\models ``\gd$ isn't strongly compact'',
contradicting our assumption that $V^{\FP}
\models ``\gd$ is strongly compact''.

We complete the proof of Lemma \ref{l4} by
showing that any measurable limit of $\K$ or
supercompact cardinal in $V^\FP$ has this feature
also in both $V$ and $V_0$.
We first verify that the forcing $\FP$ creates no
new measurable limits of $\K$. Suppose $\gk$ is a
limit point of $\K$ that
is not measurable in $V$. If $\gk$ is not a limit
point of ${\cal A}$, let $\eta$ be the least
element of ${\cal A}$ above $\gk$. As in Lemma
\ref{l3}, factor $\FP$ as $\FQ^\eta \times
\FP_\eta \times \FQ_{< \eta}$. Note that $\gk \in
(\gd_\eta, \eta)$, and as in Lemma \ref{l3},
$\FP_\eta$ satisfies the same definition in
either $V$ or $V^{\FQ^\eta}$. In particular, we
know that $\FQ^\eta$ is $\eta^+$-directed closed,
and $\FP_\eta$ admits a gap in either $V$ or
$V^{\FQ^\eta}$ below the least inaccessible above
$\eta$. Putting these facts together, and using
the Gap Forcing Theorem of \cite{H2} and
\cite{H3}, we may therefore infer that
$V^{\FQ^\eta \times \FP_\eta} \models ``\gk$
isn't measurable''. As $|\FQ_{< \eta}| < \gk$,
the results of \cite{LS} then immediately tell us
that $V^{\FQ^\eta \times \FP_\eta \times \FQ_{<
\eta}} = V^\FP \models ``\gk$ isn't measurable''
as well.

Suppose now that $\gk$ is a limit point of ${\cal
A}$. This allows us to factor the forcing as $\FP
= \FQ^\gk \times \FQ_{{<} \gk}$, where $\FQ^\gk =
\prod_{\gb
> \gk, \gb \in {\cal A}} \FP_\gb$, and
$\FQ_{< \gk} = \prod_{\gb < \gk, \gb \in {\cal
A}} \FP_\gb$. Since $\FQ^\gk$ is
${(2^\gk)}^+$-directed closed, it does not affect
whether $\gk$ is measurable, and so $\gk$ is not
measurable in $V^{\FQ^\gk}$.
Since the forcing $\FQ_{{<} \gk}$ is an Easton
support product in both $V$ and $V^{\FQ^\gk}$
(and is therefore $\gk$-c.c$.$ in both $V$ and
$V^{\FQ^\gk}$), it follows by the argument given
in the proof of Lemma 8 of \cite{A97} that $\gk$
is not measurable in $V^{\FQ^\gk \times \FQ_{{<}
\gk}} = V^\FP$.
Thus, if $\gk$ is a measurable limit point of
$\K$ in $V^\FP$, it must be one in $V$ also. As
$\FP^*$ can be factored as $\FQ_0 \ast \dot
\FQ_1$, where $|\FQ_0| = \go$ and
$\forces_{\FQ_0} ``\dot \FQ_1$ is
$\ha_1$-strategically closed'', the Gap Forcing
Theorem of \cite{H2} and \cite{H3} implies that
$\gk$ is measurable in $V_0$ as well.

Assume now that $\gk$ is supercompact in $V^\FP$.
Factor $\FP^* \ast \dot \FP$ as $\FQ_2 \ast \dot
\FQ_3$, where $|\FQ_2| = \go$ and
$\forces_{\FQ_2} ``\dot \FQ_3$ is
$\ha_1$-strategically closed''. Since $\FP^* \ast
\dot \FP$ therefore admits a gap at $\ha_1$, the
Gap Forcing Theorem of \cite{H2} and \cite{H3}
once again tells us that any supercompact
cardinal in $V^\FP = V^{\FP^* \ast \dot \FP}_0$
had to have been an element of $\K$ in $V_0$.
This, together with the remarks given in the
first paragraph of the proof of the Main Theorem
(which tell us that supercompactness is preserved
in $V$), completes the proof of Lemma \ref{l4}.
\end{proof}

This completes the proof our Main Theorem,
Theorem \ref{t1}.
\end{proof}

Let us now turn to the Generalized Main Theorem.

\begin{proof}
We will show that the appropriate conclusions of
the Main Theorem hold, more generally, for any
class $\cal A$ of supercompact cardinals having
finite Cantor-Bendixon rank. The proof will
proceed by induction on the Cantor-Bendixon rank
of $\cal A$.  If $\cal A$ has rank $1$, then it
contains none of its limit points, and the
previous theorem applies. Consider now a class
$\cal A$ having rank $n+1$. The first
Cantor-Bendixon derivative ${\cal B}={\cal
A}^{(1)}={\cal A}\intersect {\cal A}'$ has rank
$n$, and so by the induction hypothesis, there is
a forcing extension $V_{\cal B}=V^{\FP_{\cal B}}$
in which the cardinals of $\cal B$ become
non-supercompact strongly compact cardinals, and
all the supercompact cardinals of ${\cal K}-{\cal
B}$ are preserved.

Our strategy is simply to apply the Main Theorem
in the model $V_{\cal B}$ to the remaining
cardinals in $\cal A$, that is, to the cardinals
in ${\cal C}={\cal A}-{\cal B}$. Since these are
precisely the elements of $\cal A$ that are
isolated in $\cal A$, the class $\cal C$ contains
none of its limit points. Furthermore, in
$V_{\cal B}$, the class of supercompact cardinals
is precisely ${\cal K_B}={\cal K}-{\cal B}$, and
$\cal C$ is a subclass of $\cal K_B$. Thus, the
Main Theorem applies in $V_{\cal B}$ to yield a
further forcing extension $V_{\cal A}$ in which
the cardinals of $\cal C$ become non-supercompact
strongly compact cardinals and the cardinals of
${\cal K_B}-{\cal C}$ remain supercompact.
Furthermore, all the previously prepared
non-supercompact strongly compact cardinals in
$\cal B$ remain measurable limits of cardinals in
${\cal C} \subseteq {\cal A}$, since inductively,
they satisfy GCH in $V_{\cal B}$, and so they
remain strongly compact. They do not become
supercompact again because the forcing of the
Main Theorem does not create supercompact
cardinals.

In summary, the cardinals in ${\cal K}-{\cal A}$
remain supercompact through both steps of the
forcing, and the cardinals in $\cal A$ become
non-supercompact strongly compact cardinals
either in the first step, if they are limits of
$\cal A$, or in the second step, if they are
isolated in $\cal A$, respectively.
\end{proof}


The Generalized Main Corollary now follows,
because if the entire class $\cal K$ has finite
Cantor-Bendixon rank, then so also does any
subclass ${\cal A}\subseteq{\cal K}$, and so the
conclusions of the Main Theorem would hold for
any class $\cal A$.  In particular, by cutting
the universe off at some supercompact cardinal of
Cantor-Bendixon rank $n$, one obtains a model
with proper classes of supercompact cardinals of
rank below $n$. When $n>2$, for example, such a
model would have a proper class of supercompact
limits of supercompact limits of supercompact
cardinals.


We conclude this paper by noting that another
interesting generalization of the Main Theorem is
obtained when there are no supercompact limits of
supercompact cardinals and ${\cal A} = \K$. If we
then use as $\FP_\gk$, for each $\gk \in {\cal
A}$, either of the partial orderings as described
at the end of Section \ref{s2}, we obtain a
forcing extension where $V^\FP \models ``$ZFC +
There is a proper class of strongly compact
cardinals + No strongly compact cardinal $\gk$ is
supercompact + Every strongly compact cardinal
$\gk$ is $\gk^+$ supercompact and has trivial
Mitchell rank with respect to $\gk^+$
supercompactness''. This sort of model was first
constructed on pages 113--114 of \cite{A97}, but
from the much stronger hypothesis of ``ZFC + GCH
+ There is an inaccessible limit of cardinals
$\gd$ which are both $\gd^+$ supercompact and a
limit of supercompact cardinals''.

This generalizes the theorem from \cite{A95}
(which is in itself a generalization of Theorem 2
of \cite{A80}). In this result, starting from $n
\in \omega$ supercompact cardinals $\gk_1,
\ldots, \gk_n$, a model for the theory ``ZFC +
$\gk_1, \ldots, \gk_n$ are the first $n$ strongly
compact cardinals + For $1 \le i \le n$, $\gk_i$
isn't supercompact but is $\gk^+_i$
supercompact'' is constructed. Further
generalizations, e.g., producing models in which
there is a proper class of strongly compact
cardinals, no strongly compact cardinal $\gk$ is
supercompact, yet every strongly compact cardinal
$\gk$ is $\gk^{++}$ supercompact and has trivial
Mitchell rank with respect to $\gk^{++}$
supercompactness, etc., are also possible.

\end{document}